\documentclass[parskip=half,bibliography=totoc]{scrartcl}

\usepackage[automark]{scrlayer-scrpage}								%KOMA styles					
\usepackage{amsfonts}												%Mathematics fonts
\usepackage{amsmath}												%general math tools
\usepackage{mathtools}												%General mathematics symbols
\usepackage{amssymb}												%More symbols
\usepackage{extarrows}												%Extendible arrows
\usepackage{dsfont} 												%Identity matrix symbol
\usepackage{mathrsfs}												%To get mathscr
\usepackage{accents}												%Accents on math symbols
\usepackage[T1]{fontenc}											%Accents output improvement
\usepackage[utf8]{inputenc}											%Accents input improvement
\usepackage[top=1in,bottom=1in,left=1.25in,right=1.25in]{geometry}	%margins
\usepackage{enumerate} 												%Numbered lists
\usepackage{authblk}												%author/affiliation formatting
\usepackage{abstract}												%abstract
\usepackage{xcolor}	
\usepackage{tikz-cd}												%Commutative Diagrams

\usepackage{amsthm}													%Theorems etc

\usepackage[colorlinks=true,citecolor=blue]{hyperref}				%hyperlinks
\usepackage[noabbrev]{cleveref}										%better cross-referencing
\crefname{lem}{lemma}{lemmata}

\newcommand{\p}{\partial}
\newcommand{\pd}[2]{\frac{\partial #1}{\partial #2}}
\renewcommand{\d}{\mathrm{d}}
\newcommand{\dif}{d}

\newcommand{\CP}{\mathbb{C}\mathrm{P}}
\newcommand{\CH}{\mathbb{C} H}
\newcommand{\N}{\mathbb{N}}
\newcommand{\R}{\mathbb{R}}
\newcommand{\C}{\mathbb{C}}
\renewcommand{\H}{\mathbb{H}}

\renewcommand{\Im}{\operatorname{Im}}

\DeclareMathOperator{\End}{End}
\DeclareMathOperator{\id}{id}

\DeclareMathOperator{\Heis}{Heis}

\newcommand{\h}{\mathrm{H}}

\newcommand{\abs}[1]{\lvert #1 \rvert}

\newcommand{\mf}[1]{\mathfrak{#1}}
\newcommand{\mc}[1]{\mathcal{#1}}
\newcommand{\ms}[1]{\mathscr{#1}}

\newcommand{\sslash}{\mathbin{\hspace{-2pt}/\mkern-6mu/\hspace{-2pt}}}

\newtheoremstyle{mythm}% name of the style to be used
{}% measure of space to leave above the theorem. E.g.: 3pt
{}% measure of space to leave below the theorem. E.g.: 3pt
{\slshape}% name of font to use in the body of the theorem
{}% measure of space to indent
{\bfseries\sffamily}% name of head font
{.}% punctuation between head and body
{ }% space after theorem head; " " = normal interword space
{}% Manually specify head
\newtheoremstyle{mydef}% name of the style to be used
{}% measure of space to leave above the theorem. E.g.: 3pt
{}% measure of space to leave below the theorem. E.g.: 3pt
{}% name of font to use in the body of the theorem
{}% measure of space to indent
{\bfseries\sffamily}% name of head font
{.}% punctuation between head and body
{ }% space after theorem head; " " = normal interword space
{}% Manually specify head

\theoremstyle{mythm}
\newtheorem{thm}{Theorem}[section]
\newtheorem{prop}[thm]{Proposition}

\newtheorem{lem}[thm]{Lemma}
\theoremstyle{mydef}
\newtheorem{mydef}[thm]{Definition}
\newtheorem{rem}[thm]{Remark}

\newtheorem*{nota}{Notation}
\newtheorem{thm*}{Theorem}

\newenvironment{myproof}[1][\proofname]{
	\proof[\sffamily\upshape#1]
}{\endproof}

\newcommand{\proofclear}{\hfill \qedsymbol}

\clearpairofpagestyles
\ihead[]{\headmark}
\ohead[]{\pagemark}
\newenvironment{numberedlist}{\begin{enumerate}[\upshape(i)]}{\end{enumerate}}
\newenvironment{letteredlist}{\begin{enumerate}[\upshape a)]}{\end{enumerate}}

\title{Hermitian structures on a class of quaternionic K\"ahler manifolds}
\author{Vicente Cort\'es}
\author{Arpan Saha}
\author{Daniel Thung}
\affil{\normalsize  University of Hamburg\\
	vicente.cortes@uni-hamburg.de, daniel.thung@uni-hamburg.de, arpan.saha@icmat.es}
\date{}

\date{\large \today}

\begin{document}

\maketitle

\begin{abstract}
	Any quaternionic K\"ahler manifold $(\bar N,g_{\bar N},\mathcal Q)$ equipped with a Killing vector field $X$ with nowhere vanishing quaternionic moment map carries an integrable almost complex structure $J_1$ that is a section of the quaternionic structure $\mathcal Q$. Using the HK/QK correspondence, we study properties of the almost Hermitian structure $(g_{\bar N},\tilde J_1)$ obtained by changing the sign of $J_1$ on the distribution spanned by $X$ and $J_1X$. In particular, we derive necessary and sufficient conditions for its integrability and for it being conformally K\"ahler. We show that for a large class of quaternionic K\"ahler manifolds containing the one-loop deformed c-map spaces, the structure $\tilde J_1$ is integrable. We do also show that the integrability of $\tilde J_1$ implies that $(g_{\bar N},\tilde J_1)$ is conformally K\"ahler in dimension four, but not in higher dimensions.
	 In the special case of the one-loop deformation of the quaternionic K\"ahler symmetric spaces dual to the complex Grassmannians of two-planes we construct a third canonical Hermitian structure $(g_{\bar N},\hat J_1)$. Finally, we give a complete local classification of quaternionic K\"ahler four-folds for which $\tilde J_1$ is integrable and show that these are either locally symmetric or carry a cohomogeneity $1$ isometric action generated by one of the Lie algebras $\mathfrak{o}(2)\ltimes\mathfrak{heis}_3(\mathbb R)$, $\mathfrak{u}(2)$, or $\mathfrak{u}(1,1)$.
	\par
	\emph{Keywords: quaternionic K\"ahler manifolds, HK/QK correspondence, $c$-map, complex structure, conformally K\"ahler, cohomogeneity one.}\par
	\emph{MSC classification: 53C26.}
\end{abstract}

\tableofcontents

\section{Introduction}

Despite decades of effort, quaternionic K\"ahler manifolds remain arguably the most elusive manifolds of special holonomy. In particular, for a long time, only a small number of constructions of complete examples beyond locally symmetric spaces was known \cite{Ale1975,LeB1991}. Though it is conjectured that no such examples are forthcoming in the case of positive scalar curvature \cite{LS1994,BWW2020}, much progress has been made in recent years in the case of negative scalar curvature. The main driving force behind these developments has been a construction known as the (supergravity) $c$-map, and its (one-loop) deformation (see, for instance, \cite{ACDM2015} and references therein), which gives rise to many explicit examples of complete quaternionic K\"ahler manifolds of negative scalar curvature. 

The input for the $c$-map is what is known as a projective special K\"ahler (PSK) manifold, or equivalently a conical affine special (pseudo-)K\"ahler (CASK) manifold $M$. CASK manifolds have the interesting property that their cotangent bundle $N=T^*M$ naturally admits the structure of a pseudo-hyper-K\"ahler manifold and moreover comes equipped with a nowhere-vanishing rotating Killing field $Z$, which means that it fixes one of the complex structures defining the pseudo-hyper-K\"ahler structure but rotates the remaining two. 

In order to obtain a quaternionic K\"ahler manifold, one applies what is called the HK/QK correspondence to $N$ \cite{Hay2008,ACM2013,ACDM2015}. This construction, which applies to any pseudo-hyper-K\"ahler manifold with rotating Killing field $Z$ along with a choice of Hamiltonian for $Z$ with respect to the invariant complex structure, produces a pseudo-quaternionic K\"ahler manifold endowed with a Killing field $X$. The conditions 
for the resulting metric to be positive-definite are known in general \cite{ACM2013} and are satisfied in the 
case of the cotangent bundles of CASK manifolds.

Since the Killing field $X$, present on any quaternionic K\"ahler manifold $(\bar N,g_{\bar N},\mc Q)$ in the image of the HK/QK correspondence, has nowhere-vanishing quaternionic moment map \cite{Dyc2015}, $\bar N$ always carries an integrable almost complex structure \cite[\S~7]{Sal1999}. This complex structure, which we call $J_1$, is compatible with the quaternionic structure $\mc Q$. Changing the sign of $J_1$ on the two-dimensional distribution spanned by $X$ and $J_1 X$ yields another natural almost complex structure, which we call $\tilde J_1$. Note that $\tilde J_1$ is never compatible with $\mc Q$ since it induces the opposite orientation. 

In this paper, we study the properties of the almost Hermitian structure $(\tilde J_1,g_{\bar N})$. In particular, we are interested to know when  $\tilde J_1$ is also integrable, in which case we refer to the quaternionic Kähler manifold as a \emph{doubly integrable HK/QK manifold.} Accordingly,  our first main result, \Cref{thm:Hermitian}, is a necessary and sufficient condition for the integrability of $\tilde J_1$:

\begin{thm*}\label{thm:intro1}
	Let $(\bar N,g_{\bar N},\mc Q)$ be a (pseudo-)quaternionic K\"ahler manifold that arises from the HK/QK correspondence applied to a pseudo-hyper-K\"ahler manifold $(N,g_N,I_1,I_2,I_3)$ endowed with  rotating Killing field $Z$. Then the almost complex structure $\tilde J_1$ on $\bar N$ is integrable if and only if there exists a smooth function $\psi$ on the dual pseudo-hyper-K\"ahler manifold $N$ such that the identity $\nabla_Z Z=\psi I_1 Z$ holds.
\end{thm*}

In the case where $\bar N$ arises from the $c$-map, the equation $\nabla_Z Z=\psi I_1 Z$ is satisfied with $\psi\equiv -1$ (see \cite[Prop.~2]{ACM2013})  and therefore the class of doubly integrable HK/QK manifolds includes the (one-loop deformed) $c$-map spaces. This extends earlier results from \cite{CDMV2015}. 

Given a doubly integrable HK/QK manifold, one may inquire as to further properties of the Hermitian structure $(\tilde J_1,g_{\bar N})$. On the one hand, it is known that no quaternionic K\"ahler manifold, with exception of manifolds locally isometric to the Grassmannian of complex $2$-planes or to its non-compact dual, is K\"ahler. On the other hand, based on the 
study of Einstein metrics conformal to  K\"ahler metrics, see \cite{DM2003} and references therein, one may suspect that, at least in some cases, the quaternionic K\"ahler metric could be conformal to a K\"ahler metric. We prove  that this is indeed the case in dimension four (cf.~\Cref{thm:cKdim4}), but never so in higher dimensions (cf.~\Cref{thm:cKhighdim}): 

\begin{thm*}\label{thm:intro2}
 Let $(\bar N,g_{\bar N},\mc Q)$ be a  doubly integrable HK/QK manifold. Then $g_{\bar N}$ is conformal to a K\"ahler metric compatible with the complex structure $\tilde J_1$ and invariant under the canonical Killing field $X$ if and only if $\bar N$ is a four-fold.
\end{thm*}

While (deformed) $c$-map spaces are doubly integrable HK/QK manifolds, with $(\tilde J_1,g_{\bar N})$ moreover being conformally Kähler in  dimension four by virtue of the above result, they do not exhaust the class of doubly integrable HK/QK manifolds. To prove this, we show in \Cref{thm:PTcK} that every four-fold in this class is locally isometric to a member of
 a three-dimensional family $g^{a,b,c}$ of metrics. Since there is only a single $c$-map space in dimension four (carrying a one-parameter family of metrics), this gives us non-trivial examples beyond the $c$-map, which we analyze in some detail. In particular, we prove that they are either locally symmetric or are of cohomogeneity $1$. In the latter case, we give an explicit action generated by a Lie algebra of Killing fields isomorphic to either $\mathfrak{o}(2)\ltimes\mathfrak{heis}_3(\mathbb R)$, $\mathfrak{u}(2)$, or $\mathfrak{u}(1,1)$ (cf.~\Cref{prop:cohom}).

The paper is organized as follows. We recall some preliminaries regarding the HK/QK correspondence and $c$-map in \Cref{sec:prelims}, after which we prove \Cref{thm:intro1} in \Cref{sec:thm1}. In \Cref{sec:third}, we study a distinguished family of examples $(\bar N_n,g^c_n)$, $n\in \N$, $c\geq 0$, for which we have a third integrable almost complex structure. The metrics $g^c_n$ are a one-parameter deformation of the non-compact symmetric space $\frac{\mathrm{SU}(n,2)}{\mathrm{S}(\mathrm U(n)\times \mathrm U(2))}$ (which corresponds to $c=0$) through complete quaternionic K\"ahler metrics. They have been extensively studied, see e.g.~\cite{CRT2021,CST2022}. We take up the question of finding a K\"ahler metric in the same conformal class as the quaternionic K\"ahler metric in \Cref{sec:thm2and3}, proving \Cref{thm:intro2}. Finally, \Cref{sec:dim4} is dedicated to an in-depth analysis of the four-dimensional case, for which we give a local classification. We conclude with a detailed examination of the local models.

{\bfseries Acknowledgements}

This research was partially funded by the Deutsche Forschungsgemeinschaft (DFG, German Research Foundation) under Germany's Excellence Strategy -- EXC 2121 Quantum Universe  -- 390833306 and the CRC 1624. A.~S.~acknowledges support by the Spanish Ministry of Science and Innovation, through the `Severo Ochoa Programme for Centres of Excellence in R\&D' (CEX2019-000904-S),  the Proyectos de I+D+i grant PID2019-109339GA-C32, and Europa Excelencia grant EUR2020-112265, as well as support by the German Research Foundation, through the DFG Emmy Noether grant AL 1407/2-1. V.~C.~is grateful to Vestislav Apostolov and Claude LeBrun for helpful discussions about Einstein metrics which are conformally K\"ahler.

\section{Quaternionic K\"ahler manifolds from the HK/QK correspondence}
\label{sec:prelims}
	
In this preparatory section we recall the HK/QK correspondence in some detail as well as the properties of the resulting quaternionic K\"ahler manifolds that we will need later on. For more details see \cite{ACDM2015,MS2015,CST2021} and references therein.

\subsection{Twists and the HK/QK correspondence}

The starting point for the HK/QK correspondence is a pseudo-hyper-K\"ahler manifold $(N,g_N,I_k)$, $k=1,2,3$, endowed with a distinguished vector field. 

\begin{mydef}
	A vector field $Z$ on a pseudo-hyper-K\"ahler manifold $(N,g_N,I_k)$ such that $g_N(Z,Z)$ is nowhere-vanishing is called a rotating Killing field if it satisfies $L_Zg_N=0$, $L_ZI_1=0$ and $L_ZI_2=I_3$.
\end{mydef}

In order to apply the HK/QK correspondence, we ask that $(N,g_N,I_k)$ comes equipped with a rotating Killing field $Z$ which is Hamiltonian with respect to the two-form $\omega_\h=\omega_1+\d(g_N(Z,-))$, where $\omega_1$ is the K\"ahler form associated to $I_1$. Furthermore, we will require that the Hamiltonian function, which we denote by $f_\h$, is nowhere-vanishing. We remark that, under these assumptions, $Z$ is automatically Hamiltonian with respect to $\omega_1$, with Hamiltonian function $f_Z\coloneqq f_\h-g(Z,Z)$. We will also assume that $f_Z$ is nowhere-vanishing.

Under the assumption that $\omega_1$ is an integral form, the triple $(Z,\omega_\h,f_\h)$ forms twist data on $N$ in the sense of A.~Swann \cite{Swa2010}. In particular, there exists a principal circle bundle $\pi_N\colon P\to N$ with connection $\eta$ whose curvature is $\omega_\h$, and $Z$ lifts to a vector field $Z_P$ on $P$ that preserves $\eta$. Explicitly, we have
\begin{equation}\label{eq:liftZ}
	Z_P=\tilde Z+\pi_N^*f_\h X_P
\end{equation}
where $X_P$ generates the principal circle action on $P$ and $\tilde Z$ is the $\eta$-horizontal lift of $Z$. The quotient space $\bar N\coloneqq P/\langle Z_P\rangle$ is then called the twist of $N$ with respect to the twist data $(Z,\omega_\h,f_\h)$. Under the assumption that $Z$ and $Z_P$ act freely and properly, it is a smooth manifold. Note that $\bar N$ comes equipped with a distinguished vector field $X$, obtained by pushing down $X_P$.

The power of the twist construction lies in the fact that $Z$-invariant structures on $N$ can be carried over to $\bar N$. For instance, if $\varphi$ is an invariant function on $N$, then $\pi^*_N\varphi\in C^\infty(P)$ is invariant under $Z_P$ and therefore induces a well-defined function $\psi$ on $\bar N$, which we will call the twist of $\varphi$. Similarly, we can push down the $\eta$-horizontal lift $\tilde U$ of an invariant vector field $U$ on $N$ to obtain its twist $V=\dif \pi_{\bar N}(\tilde U)$ on $\bar N$. Twists of arbitrary tensor fields are then defined by demanding compatibility of twisting with tensor products and contractions. By comparing lifts to $P$, one can work out explicit formulas for the twist of a tensor field without much effort.

The following Lemmata, whose proofs can be found in \cite{Swa2010}, follow from such computations:

\begin{lem}\label{lem:exterior}
	If $\alpha$ is an invariant $p$-form on $N$ and $\beta$ denotes its twist with respect to twist data $(Z,\omega_\h,f_\h)$, then $\d \beta$ is the twist of $\d \alpha-f_\h^{-1}\omega_\h\wedge\iota_Z\alpha$.\proofclear
\end{lem}

\begin{lem}\label{lem:integrable}
	The twist of an invariant complex structure $I$ on $N$ with respect to twist data $(Z,\omega_\h,f_\h)$ is integrable if and only if $\omega_\h$ of type $(1,1)$ with respect to $I$.\proofclear
\end{lem}

Returning to the setting of the HK/QK correspondence, let $(N,g,I_k,Z)$, $k=1,2,3$, be a pseudo-hyper-K\"ahler manifold endowed with a rotating Killing field $Z$, and consider the twist data $(Z,\omega_\h,f_\h)$ as above. 

\begin{prop}[{\cite[\S 2.3]{CST2022}}]\label{prop:IH}
	Let $\nabla$ denote the Levi-Civita connection of $(N,g)$. Then the endomorphism field $I_\h=I_1+2\nabla Z$ satisfies the relation $g(I_\h-,-)=\omega_\h$, is skew-symmetric with respect to $g$, and commutes with each $I_k$. In particular, the two-form $\omega_\h$ is of type $(1,1)$ with respect to each $I_k$.\proofclear
\end{prop}

The quaternionic structure bundle $\langle I_1,I_2,I_3\rangle\subset \End(TN)$ is invariant under $Z$ and therefore induces an almost quaternionic structure $\mc Q$ on the resulting twist manifold $\bar N$ (with respect to twist data $(Z,\omega_\h,f_\h)$ as above). In fact, it is known \cite{ACDM2015,MS2015} that $(\bar N,\mc Q)$ carries a compatible quaternionic K\"ahler metric $g_{\bar N}$. This metric, however, is not the twist of the pseudo-hyper-K\"ahler metric $g_N$, but rather of its so-called elementary deformation $g_\h$ \cite{MS2015}. The elementary deformation $g_\h$ is related to $g_N$ by the formula
\begin{equation}\label{eq:ggHrelation}
	g_\h= \frac{1}{f_Z}g_N|_{(\H Z)^\perp}+\frac{f_\h}{f_Z^2}g_N|_{\H Z}
\end{equation}
where $\H Z=\langle Z,I_1 Z,I_2Z,I_3Z\rangle$ denotes the quaternionic span of $Z$. We say that the quaternionic K\"ahler manifold $(\bar N,g_{\bar N},\mc Q)$ arises from the HK/QK correspondence applied to $N$. Since the signature of a metric is preserved by twisting, the quaternionic K\"ahler metric $g_{\bar N}$ is positive definite if and only if $g_\h$ is. From \eqref{eq:ggHrelation}, we see that this may happen even if $g_N$ has signature $(4k,4)$, as long as $f_Z>0$ and $f_\h<0$. This observation is crucial in applying the HK/QK correspondence to the $c$-map \cite{ACM2013}. Since we are interested primarily in positive-definite quaternionic K\"ahler metrics, we will from now on assume that $g_\h$ is Riemannian (though our results remain valid, mutatis mutandis, in the indefinite case).

Note that the above construction depends on our choice of Hamiltonian function $f_\h$, which appears explicitly in the expression for the elementary deformation $g_\h$, whose twist is the quaternionic K\"ahler metric $g_{\bar N}$. Since there is always a one-parameter freedom in choosing a Hamiltonian function, the HK/QK correspondence in fact yields a one-parameter family of quaternionic K\"ahler metrics on $\bar N$.

We note the following property of quaternionic K\"ahler manifolds resulting from the HK/QK correspondence:

\begin{prop}\label{prop:cpxstr}
	Every quaternionic K\"ahler manifold that arises from the HK/QK correspondence admits an integrable complex structure compatible with the metric and quaternionic structure.
\end{prop}
\begin{myproof}
	This follows immediately from applying \Cref{lem:integrable} to the invariant complex structure $I_1$ on the hyper-K\"ahler manifold, using \Cref{prop:IH}, and noting that $I_1$ is compatible with the elementary deformation and quaternionic structure on the hyper-K\"ahler side of the correspondence.
\end{myproof}

\subsection{The \texorpdfstring{$c$}{c}-map}

In order to put the HK/QK correspondence to use we need a way to construct pseudo-hyper-K\"ahler manifolds endowed with a rotating Killing field. The theory of special K\"ahler manifolds provides a plentiful source of examples. 

\begin{mydef}
	An affine special K\"ahler manifold is a pseudo-K\"ahler manifold $(M,g_M,J_M)$ endowed with a flat, torsion-free and symplectic connection $\nabla^{\mathrm{SK}}$ such that $\nabla^{\mathrm{SK}} J_M$ is symmetric. An affine special K\"ahler manifold is called conical (or a CASK manifold) if it admits a vector field $\xi$ such that $\nabla^{\mathrm{SK}} \xi=\nabla^{g_M}\xi=\id_{TM}$, where $\nabla^{g_M}$ denotes the Levi-Civita connection of $g_M$. Moreover, we require that $\{\xi,J_M\xi\}$ generate a principal $\C^*$-action, and that $g_M(\xi,\xi)<0$ while $g_M$ is positive-definite on $\langle \xi,J_M\xi\rangle^\perp$.
\end{mydef}

We recall how to construct a pseudo-hyper-K\"ahler manifold $N$ with rotating Killing field from a CASK manifold $(M,g_M,J_M,\nabla^{\mathrm{SK}},\xi)$ of real dimension $2n+2$. The smooth manifold underlying $N$ is $T^*M$. Now note that the special K\"ahler connection $\nabla^{\mathrm{SK}}$ induces a splitting $T(T^*M)\cong \pi^*T^*M\oplus \pi^*TM$. With respect to this splitting, we define the tensor fields 
\begin{equation}\label{eq:cotangentHKstr}
	g_N=
	\begin{pmatrix}
		g_M & 0 \\ 0 & g^{-1}_M
	\end{pmatrix}
	\qquad  
	I_1=
	\begin{pmatrix}
		J_M & 0 \\ 0 & J_M^*
	\end{pmatrix}
	\qquad  
	I_2=
	\begin{pmatrix}
		0 & -\omega_M^{-1} \\ \omega_M & 0
	\end{pmatrix}
	\qquad 
	I_3=I_1I_2
\end{equation}
where we have omitted pullbacks throughout to simplify notation. It is well-known (see e.g.~\cite{ACD2002} and references therein) that these tensor fields define a pseudo-hyper-K\"ahler structure on $N=T^*M$ of quaternionic signature $(n,1)$. This construction is known as the rigid $c$-map. Moreover, the $\nabla^{\mathrm{SK}}$-horizontal lift of $-J_M\xi$ to $N$ endows it with a canonical rotating Killing field $Z$ \cite{ACM2013}.

This construction gives rise to many examples, since there is an abundance of CASK manifolds. Indeed, if $U\subset \C^n$ is any $\C^*$-invariant domain, there is a CASK manifold associated to any generic holomorphic function on $U$ which is homogeneous of degree $2$ under the $\C^*$-action, see \cite{ACD2002} for more detailed statements. 

It follows easily from the definitions that the vector field $-J_M\xi$ on a CASK manifold $M$ is Killing and Hamiltonian with negative Hamiltonian function $\frac{1}{2}g_M(\xi,\xi)$, so we may consider the K\"ahler quotient $M\sslash S^1=\mu^{-1}\big(-\frac{1}{2}\big)/S^1$ with respect to the circle action it generates.

\begin{mydef}
	Let $(M,g_M,J_M,\nabla^{\mathrm{SK}},\xi)$ be a CASK manifold. Then the K\"ahler quotient $\bar M =M\sslash S^1$ with respect to the Hamiltonian circle action generated by $-J_M\xi$ is called a projective special K\"ahler (PSK) manifold. 
\end{mydef}

Note that the induced K\"ahler metric on a PSK manifold is always positive-definite.

By this (extrinsic) definition, giving a PSK structure on a manifold is equivalent to giving the corresponding CASK manifold, and the two manifolds can be used interchangeably for our purposes. Thus, we may say that, starting from a PSK manifold, we may pass to the associated CASK manifold and apply the rigid $c$-map construction and the HK/QK correspondence to obtain a quaternionic K\"ahler manifold. Thus, we have obtained a method to associate quaternionic K\"ahler manifolds to PSK manifolds (cf.~\cite{ACDM2015}). This construction recovers what is known as the (supergravity) $c$-map \cite{FS1990}.

Recall that we obtain not just one, but an entire one-parameter family of quaternionic K\"ahler metrics from a single hyper-K\"ahler manifold with rotating Killing vector field by varying the corresponding Hamiltonian function. In the case of the $c$-map, there is a distinguished choice of Hamiltonian. We call the corresponding quaternionic K\"ahler metric the undeformed $c$-map metric and denote it by $g_{\bar N}^0$. The quaternionic K\"ahler metrics that arise from different shifting the Hamiltonian function by a constant $c\in \R$ are known as (one-loop) deformed $c$-map metrics \cite{RSV2006} and denoted by $g_{\bar N}^c$. It turns out that these metrics have good completeness properties only if we restrict to $c\geq 0$ \cite{CDS2017}, which we will do from now on.

\section{A Hermitian structure inducing the opposite orientation}
\label{sec:thm1}

Let $(\bar{N},g_{\bar{N}},\mc Q)$ be a quaternionic K\"ahler manifold obtained by the HK/QK correspondence from a pseudo-hyper-K\"ahler manifold $(N,g_N,I_1,I_2,I_3)$ endowed with rotating Killing field $Z$, as described in \Cref{sec:prelims}. Recall that it inherits a distinguished vector field $X$ from the generator of the principal circle action on $\pi_N\colon P\to N$. This vector field, which we can also view as the twist of the vector field $-f_\h^{-1}Z$ on $N$ (cf.~Equation \eqref{eq:liftZ}), is Killing with respect to $g_{\bar N}$. Its quaternionic moment map $\mu^X$ is nowhere-vanishing and therefore defines an integrable complex structure
\begin{equation}\label{eq:momentmap} 
	J_1 \coloneqq -\frac{\mu^X}{\| \mu^X\|} \in \Gamma (\mc Q),
\end{equation}
see \cite[Rem.~4.1.6]{Dyc2015}. This is, in fact, just another way of describing the complex structure whose existence is guaranteed by \Cref{prop:cpxstr}. We shall refer to the vector field $X$ as the canonical Killing field on $\bar N$.

Now let $\mc D$ be the rank-two distribution spanned by $X$ and $J_1X$. We define a new skew-symmetric almost complex structure $\tilde{J}_1$ by 
\begin{equation*}
	\tilde{J}_1|_\mc{D} \coloneqq -J_1|_\mc{D},\quad   
	\tilde{J}_1|_{\mc{D}^\perp} := J_1|_{\mc{D}^\perp}.
\end{equation*}
The orientation defined by  $\tilde{J}_1$ is opposite to the one defined by $\mc Q$. 

\begin{thm} \label{thm:Hermitian}
	Let $\bar{N}$ be obtained from $N$ by the HK/QK correspondence as summarized above. Then the almost complex structure $\tilde J_1$ is integrable if and only if the vector field $\nabla_Z Z$ on $N$ is a section of 
	the vector bundle $\langle Z,I_1 Z\rangle$. In this case, it defines a Hermitian structure $(g_{\bar{N}},\tilde{J}_1)$ compatible with the opposite orientation.
\end{thm} 

Since $Z$ is Killing and $g_N(Z,Z)$ is nowhere-vanishing, the above condition is equivalent to requiring that $\nabla_Z Z=\psi I_1 Z$ for some smooth function $\psi$ on $N$. Contracting with $g_N$ and again using that $Z$ is Killing, this becomes $\frac{1}{2}\d (g_N(Z,Z)) = -\psi\, \d f_Z$.

It is known that this holds when the pseudo-hyper-K\"ahler manifold $N$ with rotating Killing field $Z$ arises by applying the rigid $c$-map to a CASK manifold $M$. Indeed, in \cite[Prop.~2]{ACM2013} it is shown that, in this case, the identity $\nabla_Z Z=-I_1 Z$ holds (with our conventions). Let $\bar M$ denote the PSK manifold associated to $M$. Then the corresponding one-parameter family of quaternionic K\"ahler manifolds $(\bar N,g_{\bar N}^c)$ arise from the (deformed) $c$-map applied to $\bar M$. Therefore, the following result is a corollary of the previous theorem.

\begin{thm}\label{thm:cmap}
	Any quaternionic K\"ahler manifold $(\bar{N},g_{\bar{N}})$ obtained from the (deformed) c-map construction admits a Hermitian structure $(g_{\bar{N}},\tilde{J}_1)$ compatible with the opposite orientation. \proofclear
\end{thm} 

\begin{myproof}[Proof of \Cref{thm:Hermitian}] 
	We denote by $\mc{D}^{1,0}$ and $(\mc{D}^\perp)^{1,0}$ the $i$-eigenbundles of $\tilde{J}_1$ on $\mc{D}_\C$ and $\mc{D}^\perp_\C$, respectively. $\mc D^{0,1}$ and $(\mc D^\perp)^{0,1}$ signify the corresponding $-i$-eigenbundles. The distribution $\mc{D}^{1,0}= \langle X+iJ_1X\rangle$ is involutive for dimensional reasons. We first prove that $(\mc{D}^\perp)^{1,0}$ is also involutive. 
	
	Let $Y, Z\in \Gamma(\mc{D}^\perp)$. Then $[Y-iJ_1Y,Z-iJ_1Z]$ is of type $(1,0)$ with respect to $J_1$ by integrability of $J_1$. To prove that it is also of type $(1,0)$ with respect to $\tilde{J}_1$ it suffices to check that its projection onto $\mc D^{0,1}$ vanishes, which in turn follows from 
	\begin{equation*}
		g_{\bar{N}}(X+iJ_1X,[Y-iJ_1Y,Z-iJ_1Z])=0
	\end{equation*}
	or, equivalently, 
	\begin{equation*}
		[Y,Z]-[J_1Y,J_1Z] -J_1([J_1Y,Z]+[Y,J_1Z]) \perp \mc{D}.
	\end{equation*}
	Defining the one-forms $\alpha=g_{\bar N}(X,-)$ and $\beta=g_{\bar N}(J_1 X,-)$, this amounts to showing that $\d\alpha|_{{\bigwedge}^2 \mc{D}}$ and $\d\beta|_{{\bigwedge}^2 \mc{D}}$ are of type $(1,1)$ with respect to $J_1$. For $\alpha$, this follows by combining the fact that $\mu^X$ is proportional to $(\nabla X)^{\mf{sp}_1}$ and the definition of $J_1$ (see \eqref{eq:momentmap}) to deduce that $\nabla X$ commutes with $J_1$. To show it for $\beta$, we remark that the vanishing of the Nijenhuis tensor of $J_1$ implies
	\begin{equation*}
		J_1([X,Y]-[J_1X,J_1Y]) = [X,J_1Y]+[J_1X,Y].
	\end{equation*}
	Therefore 
	\begin{equation*}
		\beta ([Y,Z]-[J_1Y,J_1Z]) = -\alpha ([Y,J_1Z]+[J_1Y,Z]).
	\end{equation*}
	So the claimed property for $\beta$ follows from that for $\alpha$.  
	
	To establish the integrability of $\tilde J_1$ it now suffices to prove that $X+iJ_1X$ preserves $(\mc{D}^\perp)^{1,0}$. Note that $X$ and $J_1X$ are real holomorphic vector fields with respect to both $J_1$ and $\tilde{J}_1$. In addition, $X$ is Killing and therefore preserves not only $\mc{D}$ but also $\mc{D}^\perp$. As a consequence $X$ preserves the complex distribution $(\mc{D}^\perp)^{1,0}$. To show that this distribution is also $J_1X$-invariant it suffices to show that $J_1X$ preserves $\mc{D}^\perp$. This means that we have to verify that $\alpha([J_1X,Y])=\beta([J_1 X,Y])=0$ for all $Y\perp \mc{D}$. Using that $J_1X$ is $J_1$-holomorphic, one can easily see that the second condition actually follows from the first.
	
	Decomposing $\nabla X$ into its part $(\nabla X)^{\mathfrak{sp}_1}=\mu^X$ in $\mc Q$, which is proportional to $J_1$, and the part $(\nabla X)^{\mathfrak{sp}_n}$ ($n$ being the quaternionic dimension of $\bar N$) in the centralizer of $\mc Q$ we see that (for all $Y\perp \mc D$)
	\begin{equation*}
		\alpha ([J_1X,Y])= -\d\alpha (J_1X,Y)= -2g_{\bar{N}}((\nabla X)^{\mathfrak{sp}_n}J_1X,Y)
		=2g_{\bar{N}}((\nabla X)^{\mathfrak{sp}_n}X,J_1Y).
	\end{equation*}
	Thus, $\tilde J_1$ is integrable if and only if $(\nabla X)^{\mf{sp}_n}X$ is a section of $\mc{D}$ or, equivalently, if $\nabla_XX$ is. The latter is equivalent to $\d(g_{\bar{N}}(X,X))|_{(J_1X)^\perp} =0$, since $Xg_{\bar{N}}(X,X)=0$.
	\begin{lem} 
		The condition $\nabla_Z Z\in \Gamma(\langle Z,I_1 Z\rangle)$ is equivalent to $\d(g_{\bar{N}}(X,X))|_{(J_1X)^\perp} =0$.  
	\end{lem} 
	\begin{myproof} 
		Since $g(\nabla_Z Z,Z)=0$ by the Killing equation, the condition $\nabla_Z Z\in\Gamma(\langle Z,I_1 Z\rangle)$ is equivalent to $\d(g_N(Z,Z))|_{(I_1Z)^\perp} =0$. To see that the latter is in turn equivalent to $\d(g_{\bar{N}}(X,X))|_{(J_1X)^\perp}=0$ it suffices to go through the HK/QK correspondence.
		
		It follows from \eqref{eq:liftZ} that $X$ is the twist of the vector field $-f_\h^{-1}Z$ on $N$. Since $g_{\bar N}$ is the twist of $g_\h$, it follows that the function $g_{\bar{N}}(X,X)$ is the twist of $\frac{1}{f_\h^2}g_\h(Z,Z)$. Recall that $g_\h$ differs from $g_N$ by a rescaling on $\H Z$ and on $(\H Z)^\perp$, with factors involving only the $Z$-invariant functions $f_Z$ and $f_\h$ (cf.~Equation \eqref{eq:ggHrelation}). Since $\d f_Z=-\iota_Z\omega_1$ and  $\d (g_N(Z,Z))$ vanish on $(I_1Z)^\perp$, the same holds for $\d\big(\frac{1}{f_\h^2}g_\h(Z,Z)\big)$. Now, the fact that the complex structure $J_1$ is the twist of $I_1$ implies that the twist $g_{\bar{N}}(X,X)$ of the function $\frac{1}{f_\h^2} g_\h(Z,Z)$ satisfies $\d (g_{\bar{N}}(X,X))|_{(J_1X)^\perp}=0$ as well. 
	\end{myproof}
	This establishes the equivalence between the integrability of $\tilde J_1$ and the condition given in the statement of \Cref{thm:Hermitian}. 
\end{myproof}

It was already known that $\tilde J_1$ is integrable if the quaternionic K\"ahler manifold $(\bar N,g_{\bar N})$ arises from the \emph{undeformed} $c$-map \cite[Thm.~1 (a)]{CDMV2015}, so we may think of \Cref{thm:cmap} as an extension of this result to the one-loop deformed $c$-map. \Cref{thm:Hermitian} and its proof put it in the natural, broader context of the HK/QK correspondence.

There are a number of reformulations of the integrability condition, some of which will be of use later on:

\begin{lem}\label{lem:reformulations}
	The following conditions are equivalent: 
	\begin{letteredlist}
		\item $\nabla_Z Z$ is a section of the rank-two distribution $\langle Z,I_1 Z\rangle$.
		\item $\d (g_N(Z,Z))\wedge \d f_Z=0$.
		\item $\d f_\h\wedge \d f_Z=0$.
		\item There exists a smooth function $\psi$ such that one of the following equivalent identities holds:
		\begin{numberedlist}
			\item $\nabla_Z Z=\psi I_1Z$.
			\item $I_\h Z=(1+2\psi) I_1 Z$.
			\item $\d (g_N(Z,Z))=2\psi \d f_Z$.
			\item $\d f_\h=(1+2\psi)\d f_Z$.
		\end{numberedlist}
	\end{letteredlist}
\end{lem}
\begin{myproof}
	First, let us check the equivalence of the conditions listed in d). That (i) and (ii) are equivalent follows from $I_\h=I_1+2\nabla Z$. For an arbitrary vector field $Y$, the following identity holds:
	\begin{equation*}
		g_N(\nabla_Z Z,Y)=-g_N(\nabla_Y Z,Z)=-\tfrac{1}{2}\d (g_N(Z,Z))(Y).
	\end{equation*}
	 Now (i) holds if and only if $g_N(\nabla_Z Z,Y)=-\psi\,\d f_Z(Y)$, which is to say  $\d (g_N(Z,Z))=2\psi\,\d f_Z$, i.e.~condition (iii) holds. Moreover, (iii) $\Longleftrightarrow$ (iv) since $f_\h=f_Z+g_N(Z,Z)$. 
	 
	 It is obvious that b) and c) are equivalent to d), (iii) and d), (iv), respectively. Finally, the equivalence of a) and d), (i) follows directly from the Killing equation.
\end{myproof}

\section{Examples carrying a third Hermitian structure}
\label{sec:third}

Among the simplest examples of PSK manifolds are the complex hyperbolic spaces $\CH^n$. The corresponding CASK manifolds are
\begin{equation*}
	M_{n+1}\coloneqq \bigg\{(z_0,\dots,z_n)\in \C^{n+1}\bigg|\, \abs{z_0}^2>\sum_{j=1}^n \abs{z_j}^2\bigg\}
\end{equation*}
equipped with the flat K\"ahler structure defined by the Hermitian form
\begin{equation*}
	h=-\d z_0\otimes \d \bar z_0+\sum_{j=1}^n \d z_j\otimes\d \bar z_j.
\end{equation*}
The conical structure is induced by the standard $\C^*$-action by complex multiplication. In particular, we have
\begin{equation*}
	-J\xi =-i \sum_{j=0}^n \bigg(z_j \pd{}{z_j}-\bar z_j \pd{}{\bar z_j}\bigg).
\end{equation*}
Applying the $c$-map to $\CH^n$, we obtain a one-parameter family of complete quaternionic K\"ahler manifolds $(\bar N_{n+1},g_{n+1}^c)$, $c\geq 0$, of dimension $4n+4$. The undeformed $c$-map metric yields a a quaternionic K\"ahler symmetric space, while the deformed $c$-map metrics are of cohomogeneity one \cite{CST2021,CST2022}. We denote the corresponding hyper-K\"ahler manifold by $N_{n+1}$. 

\begin{prop}
	The endomorphism field $I_\h$ (defined in \Cref{prop:IH}) on $N_n$ is an integrable almost complex structure inducing the same orientation as $\mc Q$ if $n$ is even and the opposite orientation if $n$ is odd.
\end{prop} 
\begin{myproof}
	Observe that $N_n\cong M_n\times \C^n$ is nothing but an open subset of a quaternionic vector space, endowed with the standard hyper-K\"ahler structure of quaternionic signature $(n-1,1)$. In particular, $N_n$ is a product of complex manifolds and $I_1$ is the induced complex structure corresponding to this product structure. Using coordinates $(z_0,\dots,z_{n-1},w_0,\dots,w_{n-1})$ adapted to the product structure, we see that the rotating Killing field on $N_n$ is given by
	\begin{equation*}
		Z=-i \sum_{j=0}^{n-1} \bigg(z_j\pd{}{z_j}-\bar z_j \pd{}{\bar z_h}\bigg).
	\end{equation*}
	Decomposing $TN_n\cong TM_n\oplus T\C^n$, one now easily sees that 
	\begin{equation*}
		\nabla Z= 
		\begin{pmatrix}
			-I_1|_{TM_n} & 0 \\ 0 & 0 
		\end{pmatrix}
	\end{equation*}
	and correspondingly
	\begin{equation*}
		I_\h = I_1+2\nabla Z=
		\begin{pmatrix}
			- I_1|_{TM_n} & 0 \\ 0 & I_1|_{T\C^n}
		\end{pmatrix}.
	\end{equation*}
	In summary, $I_\h$ the natural induced complex structure on the product $\overline{M_n}\times \C^n$ of complex manifolds, where $\overline{M_n}$ denotes the conjugate of the complex manifold $M_n$. In particular, $I_\h$ is integrable.
	
	The claim regarding the orientation follows immediately from the fact that $TM_n\subset TN_n$ is a distribution of real dimension $2n$.
\end{myproof}

\begin{thm}\label{thm:third}
	For any $n\geq 2$, the quaternionic K\"ahler manifolds $(\bar N_n,g_n^c)$, $c\geq 0$, that arise from the $c$-map applied to $\CH^{n-1}$ admit a third complex structure $\hat J$. The associated Hermitian structure $(g_n^c,\hat J)$ induces the same orientation as $\mc Q$ if and only if $n$ is even.
\end{thm}

\begin{rem}\label{rem:n=1}
	In the case $n=1$, where one applies the $c$-map to a single point, the complex structures $\hat J$ and $\tilde J_1$ coincide. In higher dimensions they do not, since $\hat J$ disagrees with $J_1$ by a sign on a half-dimensional distribution while $\tilde J_1$ disagrees with $J_1$ on a two-dimensional distribution.
\end{rem}

\begin{myproof}[Proof of \Cref{thm:third}]
	It should not be surprising that we define $\hat J$ as the twist of the integrable complex structure $I_\h$. Since $I_\h$ is skew with respect to $g_{N_n}$ by \Cref{prop:IH}, the twisting two-form $\omega_\h$ is of type $(1,1)$ with respect to $I_\h$. By \Cref{lem:integrable}, $\hat J$ is then integrable. Furthermore, $I_\h$ is compatible with $g_\h$, so $\hat J$ is compatible with its twist $g_n^c$. The claim regarding the induced orientation follows immediately from the corresponding statement for $I_\h$, since twisting preserves algebraic properties.
\end{myproof}

Again, in the case of the undeformed $c$-map, the existence of the Hermitian structure defined by $\hat J$ had already been established in \cite[Thm.~2 (a)]{CDMV2015}. \Cref{thm:third} generalizes it to the deformed $c$-map and provides an interpretation in terms of the HK/QK correspondence.

\section{The conformal K\"ahler property}
\label{sec:thm2and3}

\begin{prop}\label{prop:cK}
	The members of the one-parameter family $g_1^c$, $c\geq 0$, of quaternionic K\"ahler metrics on the self-dual Einstein four-manifold $\bar N_1$ are all conformal to a K\"ahler metric compatible with the complex structure $\tilde J_1$. In other words, the Hermitian structures $(g_1^c,\tilde J_1)$, $c\geq 0$, are conformally K\"ahler.
\end{prop}
\begin{myproof}
	We start by considering the hyper-K\"ahler manifold $N_1$. Because we are in the four-dimensional case, we have $TN_1=\H Z$, so that in particular the elementary deformation consists only of a conformal rescaling: $g_\h=\frac{f_\h}{f_Z}g_N$. Correspondingly, the quaternionic K\"ahler metric $g_1^c$, with $c\geq 0$ arbitrary, is a conformal rescaling of the twist of $g_{N_1}$, which we will denote by $h^c$.
	
	Recall that a Hermitian structure is conformally K\"ahler if and only if the associated fundamental two-form $\sigma$ satisfies $\d \sigma=\theta\wedge\sigma$ for an exact one-form $\theta$. We will now show that this is the case for the structure $(h^c,\tilde J_1)$. It follows from \Cref{rem:n=1} that $\tilde J_1$ is the twist of $I_\h$, which satisfies $g_{N_1}(I_\h-,-)=\omega_\h$. This shows that the fundamental two-form $\sigma$ associated with the Hermitian structure $(h^c,\tilde J_1)$ is nothing but the twist of $\omega_\h$.
	
	We now use \Cref{lem:exterior} and conclude that the exterior derivative of $\sigma$ is the twist of the form 
	\begin{equation*}
		\d \omega_\h-\frac{1}{f_\h}\iota_Z\omega_\h \wedge \omega_\h
		=\d(\log f_\h)\wedge\omega_\h.
	\end{equation*}
	Now write $\psi$ for the twist of the $Z$-invariant function $\log f_\h$. Then, again by \Cref{lem:exterior}, $\d \psi$ is the twist of $\d \log f_\h$ and we conclude that 
	\begin{equation*}
		\d\sigma=\d\psi\wedge\sigma.
	\end{equation*}
	Since $h^c$ and $g^c_1$ are conformally equivalent, this proves that the Hermitian structure $(g^c_1,\tilde J_1)$ is conformally K\"ahler.
\end{myproof}

\begin{rem}
	In dimensions $4n$, $n\geq 2$, the hyper-K\"ahler metric $g_{N_n}$ is indefinite of signature $(4n-4,4)$ and therefore so is its twist $h^c$. The above proof still shows that $h^c$ is conformal to a pseudo-K\"ahler metric, but we can no longer use this to draw the same conclusion about $g_n^c$ since the relationship between $g_{N_n}$ and $g_\h$ is more complicated for $n>1$ (cf.~\Cref{eq:ggHrelation}).
\end{rem}

The quaternionic K\"ahler four-manifolds $(\bar N_1,g^c_1)$, $c\geq 0$, are examples of half conformally-flat Einstein four-manifolds that admit a principal isometric action of a three-dimensional Heisenberg group. These have recently been classified in arbitrary signature \cite{CM2021}. It turns out that such manifolds are always conformally K\"ahler or conformally para-K\"ahler, in accordance with \Cref{prop:cK}.

In light of the above, it is natural to ask whether every  doubly integrable HK/QK metric is conformal to a K\"ahler metric compatible with the integrable almost complex structure $\tilde J_1$. In the rest of this section, we shall show that this is indeed true in dimension four but false in higher dimensions. 
	
As earlier, our main tool is the twist interpretation of the HK/QK correspondence and our first task is to interpret the conformal Kähler condition on the doubly integrable HK/QK manifold in terms of data on the dual pseudo-hyper-Kähler manifold $(N,g_N,I_1,I_2,I_3)$ endowed with a rotating Killing field $Z$, which we know  satisfies $\nabla_Z Z=\psi I_1Z$ for some $\psi\in C^\infty(N)$, thanks to \Cref{thm:Hermitian}.

Denoting the fundamental two-form associated to the Hermitian structure $(g_{\bar N},\tilde J_1)$ by $\tilde\tau$, the condition that $(g_{\bar N},\tilde J_1)$ is conformally K\"ahler means that there exists a smooth function $\varphi\in C^\infty(\bar N)$ such that $\d\tilde\tau=-\d\varphi\wedge\tilde\tau$. Since $\tilde\tau$ is non-degenerate, the fact that $\tilde\tau$ is invariant under the canonical Killing field $X$ implies that the same holds for $\d\varphi$. 

This does not, in full generality, imply that $\varphi$ itself is invariant, but it is true under the weak assumption that the action generated by $X$ has at least one closed orbit. Indeed, since $L_X\d\varphi=\d (L_X\varphi)=0$, we know that $L_X\varphi$ is a constant (assuming that $M$ is connected) and therefore the restriction of $\varphi$ to any orbit is an affine function of the parameter along this orbit. In case of a closed orbit, this forces $\varphi$ to be constant, so $L_X\varphi=0$ on this orbit and thus everywhere. Since in our setting 
the canonical vector field $X$ on $\bar N$ is induced by the fundamental vector field $X_P$ of the principal $S^1$-bundle $P\to N$, all 
of its orbits are closed. Therefore $\varphi$ (and not only $\d\varphi$) is invariant under $X$. 

Using the invariant function $\varphi$, we can pass to its twist $\phi$ on the dual pseudo-hyper-K\"ahler manifold $N$. It follows from \Cref{lem:exterior} that $\d\phi$ is the twist of $\d\varphi$; denoting the twist of $\tilde\tau$ by $\tilde\sigma$, \Cref{lem:exterior} moreover shows that the conformal K\"ahler condition on $\bar N$ is equivalent to the equation 
\begin{equation}\label{eq:cKcondition0}
	\d\phi\wedge \tilde\sigma+\d\tilde\sigma=\frac{1}{f_\h}\omega_\h\wedge\iota_Z\tilde\sigma
\end{equation}
on $N$. To make further progress, it will be helpful to have explicit expressions at hand for $\tilde\sigma$ and its exterior derivative. We introduce some convenient notation for this purpose.

\begin{nota}
	Set $I_0=\id_{TN}$ and $\omega_0=g_N$, and define $\alpha_\mu=\omega_\mu(Z,-)$, where $\mu=0,1,2,3$. Introducing the orthogonal projection operator 
	\begin{equation*}
		\mc P_{\H Z}=\frac{1}{g_N(Z,Z)}\sum_{\mu=0}^3\alpha_\mu\otimes I_\mu Z\colon \mf X(N)\to \Gamma(\H Z)
	\end{equation*}
	and $\mc P_{(\H Z)^\perp}=1-\mc P_{\H Z}$, we define
	\begin{gather*}
		(\omega_1)_{\H Z}=\omega_1(\mc P_{\H Z}-,\mc P_{\H Z}-),\\
		(\omega_1)_{(\H Z)^\perp}=\omega_1(\mc P_{(\H Z)^\perp} -,\mc P_{(\H Z)^\perp}-).
	\end{gather*}
	Note that, since $I_1$ preserves $\H Z$, $\omega_1=(\omega_1)_{\H Z}+(\omega_1)_{(\H Z)^\perp}$.
\end{nota}
	
\begin{lem}\label{lem:sigmaexplicit}
	The two-form $\tilde\sigma$ and its exterior derivative are given by the following expressions:
	\begin{equation*}
		\tilde\sigma=\frac{f_\h}{f_Z^2g_N(Z,Z)}\big(-\alpha_0\wedge\alpha_1+\alpha_2\wedge\alpha_3\big)
		+\frac{1}{f_Z}(\omega_1)_{(\H Z)^\perp}
	\end{equation*}
	where $\omega_0=g$ and $\alpha_\mu=\omega_\mu(Z,-)$ for $j=0,1,2,3$, and 
	\begin{align*}
		\d \tilde\sigma&=\frac{f_\h}{f_Z^2g_N(Z,Z)}
		\Bigg(\bigg(-\frac{1+2\psi}{f_\h} +\frac{2}{f_Z} +\frac{2\psi}{f_\h-f_Z}\bigg)\alpha_2\wedge\alpha_3
		-\d\alpha_0\Bigg)\wedge\alpha_1\\
		&\quad +\frac{f_\h}{f_Z^2g_N(Z,Z)}(\alpha_2\wedge\omega_2+\alpha_3\wedge\omega_3)
		+\frac{1}{f_Z^2}\alpha_1\wedge (\omega_1)_{(\H Z)^\perp}+\frac{1}{f_Z}\d(\omega_1)_{(\H Z)^\perp}
	\end{align*}
	where $\nabla_Z Z=\psi I_1 Z$.
\end{lem}
\begin{myproof}
	We start from 
	\begin{equation*}
		\omega_1=(\omega_1)_{\H Z}+(\omega_1)_{(\H Z)^\perp} =\frac{1}{g_N(Z,Z)}\big(\alpha_0\wedge\alpha_1+\alpha_2\wedge\alpha_3\big)
		+(\omega_1)_{(\H Z)^\perp}
	\end{equation*}
	From here, we can obtain $\tilde\sigma=g_\h(\tilde I_1-,-)$ by elementary deformation and subsequently flipping the sign on the term that evaluates non-trivially on the distribution $\langle Z,I_1Z\rangle$. This immediately yields the claimed expression for $\tilde\sigma$. 
	
	Now we compute its exterior derivative.
	\begin{align*}
		\d \tilde\sigma&=\frac{f_\h}{f_Z^2g_N(Z,Z)}
		\bigg(\frac{\d f_\h}{f_\h}-\frac{2}{f_Z}\d f_Z
		-\frac{1}{g_N(Z,Z)}\d(g_N(Z,Z))\bigg)\wedge (-\alpha_0\wedge\alpha_1+\alpha_2\wedge\alpha_3)\\
		&\quad +\frac{f_\h}{f_Z^2g_N(Z,Z)}\big(-\d(\alpha_0\wedge\alpha_1)+\d(\alpha_2\wedge\alpha_3)\big)
		-\frac{1}{f_Z^2}\d f_Z\wedge (\omega_1)_{(\H Z)^\perp}+\frac{1}{f_Z}\d(\omega_1)_{(\H Z)^\perp}.
	\end{align*}
	From $\d f_Z=-\alpha_1$ and \Cref{lem:reformulations}, we know that $\d f_Z\wedge \alpha_1$ and $\d\alpha_1$ vanish while $\d f_\h=(1+2\psi)\d f_Z$ and $\d(g_N(Z,Z))=2\psi \d f_Z$. This yields
	\begin{align*}
		\d \tilde\sigma&=\frac{f_\h}{f_Z^2g_N(Z,Z)}
		\Bigg(\bigg(-\frac{1+2\psi}{f_\h} +\frac{2}{f_Z} +\frac{2\psi}{f_\h-f_Z}\bigg)\alpha_2\wedge\alpha_3
		-\d\alpha_0\Bigg)\wedge\alpha_1\\
		&\quad +\frac{f_\h}{f_Z^2g_N(Z,Z)}\d(\alpha_2\wedge\alpha_3)
		+\frac{1}{f_Z^2}\alpha_1\wedge (\omega_1)_{(\H Z)^\perp}+\frac{1}{f_Z}\d(\omega_1)_{(\H Z)^\perp}.
	\end{align*}
	The fact that $Z$ is a rotating Killing field furthermore implies $\d\alpha_2=\omega_3$ and $\d\alpha_3=-\omega_2$.	The claimed expression for $\d\tilde\sigma$ is directly obtained from these substitutions.
\end{myproof}

\begin{lem}\label{lem:conformal}
	Assume there exists a smooth function $\varphi\in C^\infty(\bar N)$, invariant under the canonical Killing field on $\bar N$, such that $(e^\varphi g_{\bar N},\tilde J_1)$ is a K\"ahler structure. Then, if $\phi$ denotes the twist of $\varphi$, there exists a smooth function $\xi$ on $N$ such that $\d \phi=\xi\, \d f_Z$.
\end{lem}
\begin{myproof}
	Recall \eqref{eq:cKcondition0}, which shows that the conformal K\"ahler condition is equivalent to the equation
	\begin{equation}\label{eq:cKcondition1}
		\d\phi\wedge \tilde\sigma+\d \tilde\sigma=f_Z^{-2}\omega_\h\wedge\d f_Z
	\end{equation}
	where we used that $\iota_Z\tilde\sigma=\frac{f_\h}{f_Z^2}\d f_Z$ by \Cref{lem:sigmaexplicit}.
	
	Contracting with $Z$, making use of the $Z$-invariance of $\phi$ and recalling that $\iota_Z\omega_\h=-\d f_\h$, we obtain
	\begin{equation*}
		\d\iota_Z\tilde\sigma=\d\bigg(\frac{f_\h}{f_Z^2}\bigg)\wedge \d f_Z=\frac{1}{f_Z^2}\big(\d f_\h-f_\h\d \phi\big)\wedge\d f_Z.
	\end{equation*}
	Since $\d f_\h\wedge \d f_Z=0$ (by the integrability of $\tilde J_1$, cf.~\Cref{lem:reformulations}) and $f_\h$ is nowhere-vanishing, this is equivalent to $\d \phi\wedge \d f_Z=0$ or in other words $\d \phi=\xi\,\d f_Z$ for some function $\xi$.
\end{myproof}

\begin{thm}\label{thm:cKdim4}
		Let $(\bar N,g_{\bar N},\mc Q)$ be a  doubly integrable HK/QK four-fold. Then $g_{\bar{N}}$ is conformal to a K\"ahler metric which is compatible with  $\tilde J_1$ and invariant under the canonical Killing field on $\bar N$.
\end{thm}

\begin{myproof}
	In the four-dimensional case, $TN=\H Z$. The assumption $\nabla_Z Z=\psi I_1 Z$, together with the fact that $I_\h$ commutes with the complex structures $I_k$, $k=1,2,3$, then implies that $I_\h=-(1+2\psi)\tilde I_1$ (this generalizes \Cref{rem:n=1}, which corresponds to the case $\psi=-1$). Contracting with $g_N$, which in the four-dimensional case equals $\frac{f_Z^2}{f_\h}g_\h$, we see that this is equivalent to
	\begin{equation}\label{eq:omegaH-sigma}
		g_N(I_\h-,-)=\omega_\h=-\frac{f_Z^2}{f_\h}(1+2\psi)\tilde \sigma.
	\end{equation}
	Now substituting \eqref{eq:omegaH-sigma} into \eqref{eq:cKcondition1}, writing $\d \phi=\xi\, \d f_Z$ (cf.~\Cref{lem:conformal}) and using the above observation, we get
	\begin{equation*}
		\d \tilde\sigma=-\big(\xi+f_\h^{-1}(1+2\psi)\big)\d f_Z\wedge \tilde\sigma.
	\end{equation*}
	We can obtain another expression for $\d \tilde\sigma$ by either restricting the result of \Cref{lem:sigmaexplicit} to dimension four or by realizing that, in dimension four, we have 
	\begin{equation*}
		\tilde\sigma=\frac{f_\h}{f_Z^2}\omega_1-2\frac{f_\h}{f_Z^2}\frac{\d f_Z\wedge\alpha_0}{f_\h-f_Z}.
	\end{equation*}
	Taking the exterior derivative and once again applying $\d f_\h\wedge \d f_Z=0$, as well as the identity $\d \alpha_0=\omega_\h-\omega_1$, this yields:
	\begin{align*}
		\d \tilde\sigma&=\frac{f_\h}{f_Z^2}\bigg(\bigg(\frac{\d f_\h}{f_\h}
		-2\frac{\d f_Z}{f_Z}\bigg)\wedge\omega_1
		+2\frac{\d f_Z\wedge (\omega_\h-\omega_1)}{f_\h-f_Z}\bigg)\\
		&=\frac{f_\h}{f_Z^2}\bigg(\bigg(\frac{\d f_\h}{f_\h}
		-2\frac{\d f_Z}{f_Z}-2\frac{\d f_Z}{f_\h-f_Z}\bigg)\wedge\omega_1
		+2\frac{\d f_Z}{f_\h-f_Z}\wedge \omega_\h\bigg)\\
		&=\bigg(\frac{1+2\psi}{f_\h}-\frac{2}{f_Z}
		-\frac{4+4\psi}{f_\h-f_Z}\bigg)\d f_Z\wedge \tilde\sigma
	\end{align*}
	where, in the last step, we used $\d f_Z\wedge \omega_1=\frac{f_Z^2}{f_\h}\d f_Z\wedge \tilde\sigma$ and \eqref{eq:omegaH-sigma}, as well as the fact that $\d f_\h=(1+2\psi)\d f_Z$ by \Cref{lem:reformulations}. 
	
	Comparing the two expressions for $\d\tilde\sigma$ yields an equation which we can solve for $\xi$:
	\begin{equation}\label{eq:xi}
		\xi = \frac{2}{f_Z}-\frac{2+4\psi}{f_\h}
		+\frac{4+4\psi}{f_\h-f_Z}.
	\end{equation}
	Integrating $\d \phi=\xi\,\d f_Z$ now yields $\phi$ and, after twisting, the conformal factor required to obtain a K\"ahler metric.
\end{myproof}

The obvious next question is whether this phenomenon persists in higher dimensions. The following result shows that this  is never the case:
\begin{thm}\label{thm:cKhighdim}
Let $(\bar N,g_{\bar N},\mc Q)$ be a  doubly integrable HK/QK manifold of dimension $\dim \bar N>4$. Then there is no function $\varphi$ invariant under the canonical Killing field on $\bar N$ such that $(e^\varphi g_{\bar N},\tilde J_1)$ is K\"ahler. 
\end{thm}
\begin{myproof} Assume that $(e^\varphi g_{\bar N},\tilde J_1)$ is K\"ahler for some function $\varphi$ invariant under the canonical Killing field on $\bar N$. Then by \Cref{lem:conformal} we know that $\d \phi=\xi \,\d f_Z$ for a smooth function $\xi$. By contracting equation \eqref{eq:cKcondition1} with a triple of vector fields, we may derive an expression for $\xi$. Let us consider any section $V\in \Gamma((\H Z)^\perp)$ with $g_N(V,V)\neq 0$ and evaluate on the triple $(I_1 Z,V,f_Z I_1 V)$. Since $(\d f_Z\wedge\tilde\beta)(I_1Z,A,B)=-g_N(Z,Z)\beta(A,B)$ for any two-form $\beta$ and any two vector fields $A$, $B$ perpendicular to $I_1Z$, we have to compute $\tilde\sigma(V,f_ZI_1V)$ and $\omega_\h(V,f_ZI_1V)$. They are given by
	\begin{align*}
		\tilde\sigma(V,f_Z I_1 V)&=g_N(V,V),\\
		\omega_\h(V,f_ZI_1 V)&=f_Z\big(\omega_1(V,I_1V)+\d\alpha_0(V,I_1 V)\big)=f_Z\big(g_N(V,V)+\d\alpha_0(V,I_1V)\big).
	\end{align*}
	The final term that needs to be computed is 
	\begin{equation*}
		\d\tilde\sigma(I_1Z,V,f_Z I_1V)=-\frac{f_\h}{f_Z}\d\alpha_0(V,I_1 V)+\frac{1}{f_Z}g_N(Z,Z)g_N(V,V)
		+\d(\omega_1)_{(\H Z)^\perp}(I_1Z,V,I_1V).
	\end{equation*}
	Applying these identities, \eqref{eq:cKcondition1} yields 
	\begin{align*}
		\xi&=\frac{1}{f_Z}\bigg(1+\frac{\d\alpha_0(V,I_1V)}{g_N(V,V)}\bigg)
		-\frac{f_\h}{f_Zg_N(Z,Z)}\frac{\d\alpha_0(V,I_1V)}{g_N(V,V)}\\
		&\quad +\frac{1}{f_Z}
		+		\frac{1}{g_N(Z,Z)}\frac{\d(\omega_1)_{(\H Z)^\perp}(I_1Z,V,I_1V)}{g_N(V,V)}\\
		&= \frac{2}{f_Z}-\frac{2}{g_N(Z,Z)}\frac{g_N(\nabla_V Z,I_1 V)}{g_N(V,V)} +\frac{1}{g_N(Z,Z)}\frac{\d(\omega_1)_{(\H Z)^\perp}(I_1Z,V,I_1V)}{g_N(V,V)},
	\end{align*}
	where we have used $\d\alpha_0=2g_N(\nabla Z,-)$ and $f_\h=g_N(Z,Z)+f_Z$.

	We can further simplify this expression by rewriting the final term. Note that, since $\iota_{I_1 Z}(\omega_1)_{(\H Z)^\perp}=0$, we have 
	\begin{align*}
		\d(\omega_1)_{(\H Z)^\perp}(I_1Z,V,I_1V)&=(I_1Z)(g_N(V,V))-g_N([I_1Z,V],V) -g_N([I_1Z,I_1V],I_1V)\\
		&=2g_N(\nabla_{I_1 Z} V,V) -g_N([I_1Z,V],V)-g_N([I_1Z,I_1V],I_1V).
	\end{align*}	
	Now, we use torsion-freeness and the fact that $\nabla I_1=0$ to find
	\begin{align*}
		&\d(\omega_1)_{(\H Z)^\perp}(I_1Z,V,I_1V)\\
		&=2g_N(\nabla_{I_1 Z} V,V)-g_N(\nabla_{I_1Z}V-I_1\nabla_VZ,V)-g_N(I_1\nabla_{I_1Z}V-I_1\nabla_{I_1V}Z,I_1V)\\
		&=g_N(I_1\nabla_V Z,V)+g_N(\nabla_{I_1 V} Z,V)\\
		&=2g_N(I_1\nabla_V Z,V),
	\end{align*}
	where the final step follows from the Killing equation for $Z$.
	
	In conclusion, we have shown that 
	\begin{equation*}
		\xi= 
		\frac{2}{f_Z}+\frac{4}{g_N(Z,Z)}\frac{g_N(I_1\nabla_V Z,V)}{g_N(V,V)}.
	\end{equation*}
	The fact that $\xi$ does not depend on our choice of $V$ now imposes the existence of a smooth function $\lambda$ such that 
	\begin{equation*}
		g_N(I_1 \nabla_V  Z,V)=-\lambda\cdot g_N(V,V).
	\end{equation*}
	Next, observe that the endomorphism field $I_1\circ \nabla Z$ is symmetric with respect to $g_N$. Indeed, both $I_1$ and $\nabla Z$ are skew-symmetric with respect to $g_N$--for the latter, this is the Killing equation---and furthermore commute since $[I_1,\nabla Z]=\frac{1}{2}[I_1,I_\h-I_1]=0$ by \Cref{prop:IH}. It now follows from the polarisation identity that 
	\begin{equation*}
		g_N(I_1\nabla_V Z,W)=-\lambda \cdot g_N(V,W)
	\end{equation*}
	for any other $W\in \Gamma((\H Z)^\perp)$. This shows that $I_1\circ \nabla Z|_{(\H Z)^\perp}=-\lambda \id_{(\H Z)^\perp}$ or, equivalently, $\nabla Z|_{(\H Z)^\perp}=\lambda\, I_1|_{(\H Z)^\perp}$. Since $I_\h=I_1+2\nabla Z$, we deduce that 
	\begin{equation*}
		I_\h|_{(\H Z)^\perp}=(1+2\lambda)I_1|_{(\H Z)^\perp}.
	\end{equation*}
	To avoid running afoul of \Cref{prop:IH}, which guarantees that $I_\h$ commutes with $I_2$ and $I_3$, we must have $\lambda\equiv -\frac{1}{2}$, so in particular $$\xi=\frac{2}{f_Z}+\frac{2}{g_N(Z,Z)}=\frac{2}{f_Z}+\frac{2}{f_\h - f_Z}$$ and $I_\h|_{(\H Z)^\perp} = 0$. Meanwhile since $I_\h$ commutes with all the $I_i$ and 
	we know its action on $Z$, the restriction $I_\h|_{\H Z}$ must be $-(1+2\psi)\tilde I_1|_{_{\H Z}}$. This gives us
		\begin{equation}\label{eq:omegaH-psi}
			 \omega_\h = -\frac{1+2\psi}{g_N(Z,Z)}(-\alpha_0\wedge\alpha_1+\alpha_2\wedge\alpha_3).
			 \end{equation}
	Now, evaluating the exterior derivative of $\omega_\h$ on $I_2 Z, V, I_2V$, where $V$ is, as earlier, a section $V\in \Gamma((\H Z)^\perp)$ with $g_N(V,V)\neq 0$, we get
	\begin{align*}
		\d\omega_\h(I_2 Z, V, I_2V) &= \frac{1+2\psi}{g_N(Z,Z)}\alpha_2(I_2 Z)\d\alpha_3(V, I_2V)\\
		&= -(1+2\psi)\omega_2(V, I_2V)=-(1+2\psi)g_N(V,V).
	\end{align*}
	Here we have used that $\d\alpha_3 =-\omega_2$, as in the proof of Lemma \ref{lem:sigmaexplicit}.
	Since $\omega_\h$ is closed, the above must vanish, giving us $1+2\psi=0$. So, $\omega_\h$ and $I_\h$ are zero, implying $\nabla Z = -\frac12 I_1$ and that $f_\h$ is constant (since $\d f_\h = -\iota_Z\omega_\h = 0$). Thus, we may explicitly solve $\d\phi = \xi\,\d f_Z$ to obtain $$e^\phi = \frac{Cf_Z^2}{(f_\h-f_Z)^2},$$
	for some positive constant $C$. Meanwhile \eqref{eq:cKcondition1} just becomes the statement that $e^\phi\tilde\sigma$ is closed. We shall now verify using our explicit expressions for $e^\phi$ and $\tilde \sigma$ that this is not the case.
	
	First of all, we write $\tilde\sigma$ as
	\begin{align*}
		\tilde\sigma&=\frac{1}{f_Z}\bigg( \frac{f_\h}{f_Z(f_\h - f_Z)}(-\alpha_0\wedge\alpha_1 +\alpha_2\wedge \alpha_3) +(\omega_1)_{(\H Z)^\perp}\bigg)\\
		&=\frac{1}{f^2_Z}\bigg( -\frac{f_\h+f_Z}{f_\h - f_Z}\,\alpha_0\wedge\alpha_1 +\alpha_2\wedge \alpha_3\bigg) +\frac{1}{f_Z}\omega_1.
	\end{align*}
	We see that  evaluating the exterior derivative of $e^{\phi}\tilde\sigma$ on $I_2 Z, V, I_2V$ gives
\begin{align*}
\d(e^{\phi}\tilde\sigma)(I_2 Z, V, I_2V) &= -\frac{C}{(f_\h - f_Z)^2}\alpha_2(I_2 Z)\d\alpha_3(V, I_2V)\\
&= \frac{C\omega_2(V, I_2V)}{f_\h - f_Z}=\frac{C g_N(V,V)}{f_\h - f_Z}\neq 0.
\end{align*}
Thus $\d(e^{\phi}\tilde\sigma)\neq 0$ and our result follows.
\end{myproof}

\section{Four-dimensional examples beyond the \texorpdfstring{$c$}{c}-map}
\label{sec:dim4}

\subsection{Local classification from the Przanowski--Tod ansatz}

The one-parameter family $(\bar N_1,g_1^c)$ is known as the (one-loop) deformed universal hypermultiplet in the physics literature and its members are the only four-dimensional quaternionic K\"ahler manifolds that arise from the $c$-map (namely, by applying it to a point). Mathematically, $(\bar N_1,g_1^c)$ is a one-parameter deformation of the Hermitian symmetric space $\CH^2\cong (\bar N_1,g_1^0)$. In this section, we construct further examples belonging to the class described in \Cref{thm:Hermitian}. In order to do so, we widen our scope to arbitrary quaternionic K\"ahler four-manifolds endowed with a Killing field.

Via HK/QK correspondence, any quaternionic K\"ahler manifold of dimension four equip\-ped with a Killing field (locally) arises from a hyper-K\"ahler four-manifold $(N,g_N)$ with rotating Killing field. Such a hyper-K\"ahler metric is locally homothetic to a metric satisfying the Boyer--Finley ansatz \cite{BF1982}. With respect to local coordinates $(\rho,x,y,t)$, this means that
\begin{equation*}
	g_N= K(\p_\rho u)(\d \rho^2 + 2 e^u(\d x^2 +\d y^2))
		 +\frac{4K}{\p_\rho u}
		 \bigg(\d t-\frac{1}{2}(\p_y u\,\d x-\p_x u\,\d y)\bigg)^2
\end{equation*}
with $u$ a smooth function of $\rho$, $x$, $y$ (but not $t$) satisfying the continuous Toda equation
\begin{equation*}
	(\p_x^2+\p_y^2)u =- 2\p_\rho^2(e^u).
\end{equation*}
The rotating Killing vector field is $\p_t$ and the invariant K\"ahler form is
\begin{equation*}
		\omega_1 =2K\d\rho\wedge
		\bigg(\d t-\frac{1}{2}(\p_y u\,\d x-\p_x u\,\d y)\bigg)
		2Ke^u\p_u\rho\,\d x\wedge \d y.
\end{equation*}
The Killing vector field $\p_t$ is Hamiltonian. Without loss of generality, the Hamiltonian function can be taken to be $f_Z = 2K\rho$ (shifts in $\rho$ can be adsorbed into a redefinition of the function $u$).

The result of applying the HK/QK correspondence to the Boyer--Finley ansatz is the Przanowski--Tod ansatz \cite{Prz1991}. Up to constant scaling, any four-dimensional quaternionic K\"ahler metric admitting a Killing field is locally given by
\begin{equation*}
	g_{\bar{N}} = \frac{1}{4\rho^2}\left(P\,\d \rho^2 + 2P e^u(\d x^2 +\d y^2) + \frac{1}{P}\,(\d t + \Theta)^2\right)
\end{equation*}
where $P$ is a $\p_t$-invariant smooth function and $\Theta$ is a $\p_t$-invariant $1$-form satisfying
\begin{equation*}
	P=K(\rho\,\p_\rho u -2)> 0, \quad 
	\d \Theta =\big(\p_y P\,\d x - \p_{x} P\,\d y\big)\wedge \d \rho 
	- 2\p_\rho (Pe^u)\,\d x\wedge \d y.
\end{equation*}
We have normalized $g_{\bar N}$ so that the reduced scalar curvature $\nu$ equals $\frac{2}{K}$. In particular, since $P$ is required to be positive, the function $\rho\,\p_\rho u -2$ must have a definite sign over the domain of definition. This sign is furthermore equal to the sign of $K$ and hence the sign of the scalar curvature.

Note that if we let $\zeta = x + i y$ and $f$ be a holomorphic function of $\zeta$ whose derivative is non-vanishing in the domain of definition, then the replacement
\begin{equation}\label{eq:Toda-gauge}
	u(\rho,\zeta, \overline \zeta) \mapsto u(\rho,f(\zeta),\overline{f(\zeta)}) - \ln\left\lvert\frac{\d f}{\d\zeta}\right\rvert^2
\end{equation}
doesn't change the isometry class of the metrics $g_N$ and $g_{\bar N}$. Such a replacement can therefore be regarded as a gauge transformation.

Now we can explicitly which of such quaternionic K\"ahler metrics are conformal to a K\"ahler metric compatible with $\tilde J_1$ in terms of the function $u$ up to gauge transformations. 

\begin{thm}\label{thm:PTcK}
	A quaternionic K\"ahler four-manifold which arises from the HK/QK correspondence has the structure $\tilde{J}_1$ integrable if and only if it is locally isometric to a metric of the form
	\begin{equation}\label{eq:Liouville-QK}
		\begin{split}
			g^{a,b,c} &= -\frac{K}{4\rho^2}\left(\frac{b\rho + 2c}{a\rho^2 + b\rho + c}\,\d \rho^2 + \frac{2(b\rho+2c)|\d\zeta|^2}{(1 + \frac{a}{2}|\zeta|^2)^2}\right.\\
			&\quad+\left. \frac{a\rho^2 + b\rho + c}{b\rho + 2c}\left(-\frac{\d t}{K} + \frac{b\,\mathrm{Im}(\zeta \,\d\overline\zeta)}{1 + \frac{a}{2}|\zeta|^2}\right)^2\right)
		\end{split}
	\end{equation}
	where $a,b,c$ are real constants, and $\rho$ and $\zeta$ satisfy
	\begin{equation}
		\rho > 0,\qquad a\rho^2 + b\rho + c > 0, \qquad b\rho + 2c \neq 0,\qquad 1+ \frac{a}{2}|\zeta|^2 > 0
	\end{equation}
	and $K$ is a non-zero constant for which the metric is positive-definite. These metrics arise by applying the Przanowski--Tod ansatz to
	\begin{equation}\label{eq:u}
		e^u = \frac{a\rho^2 + b\rho + c}{(1 + \frac{a}{2}|\zeta|^2)^2}
	\end{equation}
	up to gauge transformation.
\end{thm}

\begin{myproof}
	Any quaternionic K\"ahler four-manifold which arises from the HK/QK correspondence is locally homothetic to the Przanowski--Tod ansatz, so we restrict to this case. By \Cref{thm:cKdim4} and \Cref{lem:reformulations}, it suffices to show that $\d (g_N(Z,Z))\wedge \d f_Z=0$ if and only if \eqref{eq:u} holds. The Hamiltonian $f_Z$ is $4K\rho$ in our case while $g_N(Z,Z)$ can be seen to be $4K/\p_\rho u$. So, the condition above amounts to saying that the differential of $\p_\rho u$ is proportional to $\d \rho$, which is to say $u$ is a separable function
	\begin{equation*}
		u = F(\rho) + G(\zeta,\overline\zeta).
	\end{equation*}
	Substituting this into the continuous Toda equation, we obtain
	\begin{equation}
		\p_\zeta\p_{\overline\zeta} G = -ae^G, \quad \p_\rho^2 e^F = 2a.
	\end{equation}
	The equation for $F$ implies that $e^F$ is a quadratic polynomial in $\rho$ with leading coefficient $a$. That is, its general solutions are given by
	\begin{equation}\label{eq:expF}
		e^F = a\rho^2 + b\rho + c.
	\end{equation}
	Note in particular that the function $\rho$ is nowhere vanishing and so has a definite sign if we assume connectedness. Since the sign can always be absorbed into the constant $b$, we can assume without loss of generality $\rho >0$.

	Meanwhile, the equation for $G$ is just the 2d Liouville equation whose general solutions are known to be of the form
	\begin{equation}
		e^G = \frac{4}{(1 + 2a|f(\zeta)|^2)^2}\left\lvert\frac{\d f}{\d \zeta}\right\rvert^2
	\end{equation}
	where $f$ is some holomorphic function which is nonvanishing in the domain of definition. The freedom to choose $f$ may be absorbed into the gauge transformation  \eqref{eq:Toda-gauge}. In particular, we may set $f(\zeta) = \frac{1}{2}\zeta$. This gives us the required metric. 		
	 Moreover, we may additionally assume without any loss of generality that  $1+ \frac{a}{2}|\zeta|^2>0$. To see this, note that when $a\ge 0$, this is automatic, and when $a<0$, the replacement $\zeta\mapsto \frac{2}{a\zeta}$ leaves $e^G$ (and hence the metric) unchanged but swaps the regions $\pm(1+ \frac{a}{2}|\zeta|^2)>0$.
\end{myproof}

\begin{rem}
	As explained earlier, the sign of the scalar curvature of a metric satisfying the Przanowski--Tod ansatz is equal to the sign of the (definite) function $\rho\p_\rho u-2$, which in the case at hand satisfies $\rho \p_\rho u - 2 =-\frac{b\rho  +2c}{a\rho^2 + b\rho + c}$. The denominator, being equal to $e^F$ in \eqref{eq:expF}, is positive, so the sign of the scalar curvature of $g^{a,b,c}$ equals the sign of the function $-(b\rho  +2c)$.
\end{rem}

\subsection{Case-by-case analysis}\label{ssec:familiar}

The metrics $g^{a,b,c}$ have previously been studied by Ketov \cite{Ket2001}. There are many ranges for the parameters in which we recover well-known examples of quaternionic K\"ahler four-manifolds, as can be seen from the following case-by-case analysis of the metrics.

\begin{numberedlist}
	\item When $a=b=0$ we get  the real hyperbolic 4-space $\mathrm{SO}(4,1)/\mathrm{O}(4)$ presented as 
	solvable Lie group $\R_{>0}\ltimes \R^3$ (the Iwasawa subgroup of $\mathrm{SO}_0(4,1)$) with a left-invariant metric.
	\item For $a=0$, $b\neq 0$, the metrics $g^{a,b,c}$ are isometric to the (one-loop) deformed universal hypermultiplet metrics, i.e.~yield a one-parameter deformation of the complex hyperbolic plane $\mathrm{SU}(2,1)/\mathrm{U}(2)$ ($c=0$) through quaternionic K\"ahler metrics. For $c\neq 0$, the full isometry group of the deformed universal hypermultiplet is known to be $\mathrm{O}(2)\ltimes \Heis_3$, acting with cohomogeneity one \cite{CST2021}.	
	\item\label{item:abpos} For $a,b>0$, we may perform the following change of coordinates:
	\begin{equation*}
		\begin{split}
			\rho = \frac{b}{2a}\left(\frac{1}{\varrho^2}-1\right), \quad \zeta = \sqrt{\frac{2}{a}}\,\xi,\quad t = \frac{b K}{a}\,\theta
		\end{split}
	\end{equation*}
	where $0<\varrho<1$. In these coordinates, after making the identification $\theta\sim \theta + 2\pi$ and smoothly extending to $\zeta =\infty$ and $\varrho=0$, the metric $g^{a,b,c}$ becomes the Pedersen metric \cite{Ped1986} on the open unit ball in $\R^4$:
	\begin{equation*}
		\begin{split}
			g^{a,b,c}=\frac{-1}{\nu(1-\varrho^2)^2} \left(\frac{1 + k \varrho^2}{1 + k\varrho^4}\,\d\varrho^2 + \varrho^2(1+ k \varrho^2)(\varsigma^2_1 + \varsigma^2_2) + \frac{\varrho^2(1 + k \varrho^4)}{1 + k\varrho^2}\varsigma^2_3\right).
		\end{split}
	\end{equation*}
	Here, $k$ is given by $k = \frac{4ac}{b^2} -1$. The Pedersen metric is known to be invariant under the standard action of $\mathrm{U}(2)$ on $S^3$, as confirmed by the appearance of the $\mathrm{SU}(2)$-invariant $1$-forms $\varsigma_1,\varsigma_2,\varsigma_3$, which are given by
	\begin{equation*}
		\begin{split}
			\varsigma_1 &=\frac{\mathrm{Re}(e^{i \theta}\d\xi)}{1 + |\xi|^2}, \quad \varsigma_2 = \frac{\mathrm{Im}(e^{i \theta}\d\xi)}{1 + |\xi|^2}, \quad \varsigma_3 = \frac{1}{4}\,\d \theta + \frac{\Im(\xi\,\d \overline\xi)}{1 + |\xi|^2}
		\end{split}
	\end{equation*}
	with respect to the Hopf parametrisation. The combination in which they appear in the metric ensures that it is even $\mathrm{U}(2)$-invariant.
	\item \label{item:limiting} As a limiting case of the above, we have  $a >0$ and $b=0$.  If $c>0$, we carry out a change of coordinate $\rho = \frac{1}{\rho'}\sqrt{\frac{c}{a}}$ to obtain the following metric on $\R\times \R_{>0}\times S^2 \cong \R\times (\R^3\setminus \{0\}) \subset \R^4$:
\begin{equation*}
		\begin{split}
			g^{a,b,c} &= -\frac{1}{\nu}\left(\frac{\d \rho'^2}{1+\rho'^2} + (1+\rho'^2)\left(\frac{\d t}{2K}\right)^2 + \rho'^2 g_{S^2} \right)
		\end{split}
	\end{equation*}
where $g_{S^2}$ is the standard metric on a unit $S^2$. By including a line parametrised by $t$ at $\rho'=0$, this can be extended to the real hyperbolic 4-space metric on $\R^4$. 

On the other hand, if  $c>0$, we carry out a change of coordinate $\rho = \frac{1}{\rho'}\sqrt{\frac{c}{a}}$ to obtain the following metric on $\R\times (0,1)\times S^2 \subset \R^4$:
\begin{equation*}
	\begin{split}
		g^{a,b,c} &= -\frac{1}{\nu}\left(\frac{\d \rho'^2}{1-\rho'^2} + (1-\rho'^2)\left(\frac{\d t}{2K}\right)^2 + \rho'^2 g_{S^2} \right).
	\end{split}
\end{equation*}
Upon making the identification $t\sim t + 4\pi K$, the first two terms give the standard $S^2$ metric restricted to a hemisphere minus the pole, with $\sin^{-1}(\rho')$ and $\frac{t}{2K}$ being the latitude and longitude respectively. With this identification, the full metric $g^{a,b,c}$ can be extended to $S^4$ to give the standard $\mathrm O(5)$-invariant metric on it.

	\item When $b < 0< a$ and $2a\rho + b < 0$, the same change of coordinates as in case (\ref{item:abpos}) can be carried out, but now we have $\varrho > 1$, as can be seen from
	\begin{equation*}
		\varrho^2 =\frac{b}{2a\rho + b}= \frac{|b|}{|b| - 2a\rho} > 1.
	\end{equation*}
	We are therefore dealing with a metric on the complement of the closed unit ball in $\R^4$. The scalar curvature remains negative since
	\begin{equation*}
		-(b\rho  +2c) = (2a\rho + b)\rho - 2(a\rho^2 + b\rho + c) <0.
	\end{equation*}
	This does not change the fact that the metric is $\mathrm{U}(2)$-invariant, since it is given by the same formal expression as the Pedersen metric.
	\item When $b < 0 < a$ and $2a\rho + b > 0$, we instead set
	\begin{equation*}
		\begin{split}
			\rho = -\frac{b}{2a}\left(\frac{1}{\varrho^2}+1\right), \quad \zeta = \sqrt{\frac{2}{a}}\,\xi,\quad t = \frac{bK}{a}\,\theta.
		\end{split}
	\end{equation*}
	The metric now becomes
	\begin{equation*}
		\begin{split}
			g^{a,b,c}=\frac{1}{\nu(1+\varrho^2)^2} \left(\frac{1 - k \varrho^2}{1 + k\varrho^4}\,\d\varrho^2 + \varrho^2(1- k \varrho^2)(\varsigma^2_1 + \varsigma^2_2) + \frac{\varrho^2(1 + k \varrho^4)}{1 - k\varrho^2}\varsigma^2_3\right).
		\end{split}
	\end{equation*}
	For $k\le0$, the scalar curvature is positive since 
	\begin{equation*}
		-(b\rho+2c)= \frac{-b}{2a}(2a\rho + b + kb) >0.
	\end{equation*}
	These metrics are $\mathrm{U}(2)$-invariant in general (for the same reason as the Pedersen metric), but for particular values of $k$ we can say more. Setting $k=0$ gives the standard $\mathrm{O}(5)$-invariant metric on the complement of  a point in $S^4$. If we set $k = -1$ instead (which corresponds to $c= 0$) we obtain 
	\begin{equation*}
		\begin{split}
			g^{a,b,c}=\frac{1}{\nu(1+\varrho^2)^2} \left(\frac{\d\varrho^2}{1 - \varrho^2} + \varrho^2(1+ \varrho^2)(\varsigma^2_1 + \varsigma^2_2) + \varrho^2(1 - \varrho^2)\varsigma^2_3\right).
		\end{split}
	\end{equation*}
	This is the $\mathrm{U}(3)$-invariant Fubini--Study metric on the complement of the point $[z_0:z_1:z_2]=[1:0:0]$ in the complex projective plane with
	\begin{equation*}
		\frac{1}{\varrho^{2}} = \frac{2}{z_1^2+z_2^2} + 1.
	\end{equation*}
	\item\label{item:abneg} For $a,b<0$, we may perform the following change of coordinates:
	\begin{equation*}
		\begin{split}
			\rho = \frac{b}{2a}\left(\frac{1}{\varrho^2}-1\right), \quad \zeta = \sqrt{-\frac{2}{a}}\,\xi,\quad t = \frac{bK}{a}\,\theta
		\end{split}
	\end{equation*}
	where $0<\varrho<1$. In these coordinates, after making the identification $\theta\sim \theta + 2\pi$, the metric $g^{a,b,c}$ becomes 
	\begin{equation*}
		\begin{split}
			g^{a,b,c}=\frac{-1}{\nu(1-\varrho^2)^2} \left(\frac{1 + k \varrho^2}{1 + k\varrho^4}\,\d\varrho^2 - \varrho^2(1+ k \varrho^2)(\varsigma'^2_1 + \varsigma'^2_2) + \frac{\varrho^2(1 + k \varrho^4)}{1 + k\varrho^2}\varsigma'^2_3\right).
		\end{split}
	\end{equation*}
	This is formally very similar to the Pedersen metric, except that we now have $\mathrm{SU}(1,1)$-invariant $1$-forms $\varsigma'_1,\varsigma'_2,\varsigma'_3$ on $B^2\times S^1$ given by
	\begin{equation*}
		\begin{split}
			\varsigma'_1 &=\frac{\mathrm{Re}(e^{i \theta}\d\xi)}{1 - |\xi|^2}, \quad \varsigma'_2 = \frac{\mathrm{Im}(e^{i \theta}\d\xi)}{1 - |\xi|^2}, \quad \varsigma'_3 = \frac{1}{4}\,\d \theta - \frac{\Im(\xi\,\d \overline\xi)}{1 - |\xi|^2}.
		\end{split}
	\end{equation*}
	Correspondingly, the metric is invariant under the transitive action of $\mathrm{U}(1,1)$ on $B^2\times S^1$. Note that $k$ is necessarily negative in order for the metric to be definite. This can also be seen from the fact that, were $k$ non-negative, the discriminant of $a\rho^2 + b\rho + c$ would be positive. Since $a <0$, this would then imply that  $a\rho^2 + b\rho + c$ is nowhere positive.
\item As a limiting case of the above, we have $a<0$ and $b=0$. Noting that $c>0$ necessarily in order for $a\rho^2 + c$ to be positive, and carrying out the change of variable  $\rho = \frac{1}{\rho'}\sqrt{-\frac{c}{a}}$, we get
\begin{equation*}
	\begin{split}
		g^{a,b,c} &= -\frac{1}{\nu}\left((1-\rho'^2)\left(\frac{\d t}{2K}\right)^2+\frac{\d \rho'^2}{1-\rho'^2} + \frac{4 \rho'^2|\d\xi|^2}{(1 - |\xi|^2)^2}\right).
	\end{split}
\end{equation*}
Note that as in case \eqref{item:limiting} with $c<0$, upon making the identification $t\sim t + 4\pi K$, the first two terms give the standard $S^2$ metric restricted to a hemisphere minus the pole, with $\sin^{-1}(\rho')$ and $\frac{t}{2K}$ being the latitude and longitude respectively, while 
the last term is the disc form of the real hyperbolic plane  metric. The full metric $g^{a,b,c}$ is in fact the real hyperbolic $4$-space, presented in terms of the fibration over $S^2$ given by sending a point $(x_0,x_1,x_2,x_3,x_4)$ on the hyperboloid $x_0^2 -\sum_{i=1}^4x_i^2 = (-\nu)^{-1/2}$ sitting inside $\R^{1,4}$ to $(\sqrt{x_0^2 - x_1^2 -x_2^2},x_3,x_4)\in S^2\subset \R^3$.
\item When $a < 0< b$ and $2a\rho + b > 0$, the same change of coordinates as in case (\ref{item:abneg}) can be done, but now we have $\varrho > 1$, as can be seen from
\begin{equation*}
	\varrho^2 =\frac{b}{2a\rho + b}= \frac{b}{b - 2|a|\rho} > 1.
\end{equation*}
Because the metric is formally identical to the metric from case (\ref{item:abneg}), it is also $\mathrm{U}(1,1)$-invariant.
\item In the final case, $a < 0< b$ and $2a\rho + b < 0$, we instead set
\begin{equation*}
	\begin{split}
		\rho = -\frac{b}{2a}\left(\frac{1}{\varrho^2}+1\right), \quad \zeta = \sqrt{-\frac{2}{a}}\,\xi,\quad t = \frac{bK}{a}\,\theta.
	\end{split}
\end{equation*}
The metric now becomes
\begin{equation*}
	\begin{split}
		g^{a,b,c}=\frac{1}{\nu(1+\varrho^2)^2} \left(\frac{1 - k \varrho^2}{1 + k\varrho^4}\,\d\varrho^2 - \varrho^2(1- k \varrho^2)(\varsigma'^2_1 + \varsigma'^2_2) + \frac{\varrho^2(1 + k \varrho^4)}{1 - k\varrho^2}\varsigma'^2_3\right).
	\end{split}
\end{equation*}
which is once again $\mathrm{U}(1,1)$-invariant because of its general form.
\end{numberedlist}

A careful inspection of our list shows that, in every case, we have either identified the metric $g^{a,b,c}$ as locally symmetric or given an explicit action of a group of isometries whose orbits are three-dimensional. This means, in the latter case, that the metrics are of cohomogeneity at most one. 

It is known that the (one-loop) deformed universal hypermultiplet is of cohomogeneity (exactly) $1$ for non-zero values of the deformation parameter. It is therefore natural to ask whether this holds for the metrics $g^{a,b,c}$ which are not locally symmetric. We answer this question affirmatively:

\begin{prop}\label{pr:cohom1}
	The metric $g^{a,b,c}$ of \eqref{eq:Liouville-QK}  is locally symmetric when $bc(b^2-4ac)= 0$ and of cohomogeneity exactly $1$ otherwise. In the cohomogeneity $1$ case, there is a transitive locally isometric action of
	$\mf{o}(2)\ltimes\mf{heis}_3(\R)$ if $a=0$, $\mf{u}(2)$ if $a>0$, and $\mathfrak{u}(1,1)$ if $a<0$
	on the constant $\rho$ hypersurfaces.
\end{prop}
\begin{myproof}
	The condition $bc(b^2-4ac)= 0$ implies that either we have $b = 0$ or  $c=0$ or $b^2-4ac=0$. In the last case, if we assume that $b\neq 0$, this is the same as saying $k:=\frac{4ac}{b^2} -1=0$. From the case-by-case analysis in the previous subsection, we see that these yield locally symmetric spaces. 
	
	To see that the cohomogeneity of the metric is at least $1$ when $bc(b^2-4ac)\neq 0$, we note that its curvature norm 
	\begin{equation}\label{eq:curv-norm}
		\mathrm{tr}(\ms R^2) = 6\nu^2\left(1+b^2(b^2-4ac)^2\left(\frac{\rho}{b\rho + 2c}\right)^6\right)
	\end{equation}
	is an injective function of $\rho >0$ whenever  $b(b^2-4ac)\neq 0$  and $c\neq 0$.
	
	Though we have already given explicit actions of groups of isometries whose orbits are of codimension $1$ in every case, we can provide a direct proof that such a group must exist. When $b(b^2-4ac)\neq 0$, in the notation of \Cref{thm:PTcK}, the level sets of $\rho$ admit four Killing fields which span a three-dimensional distribution (at every point). Their explicit expressions are the following:
	\begin{equation*}
		\begin{split}
			&\mathrm{Re}\bigg(\partial_\zeta + \frac{a}{2}\,\overline{\zeta}^2\partial_{\overline\zeta} - \frac{ i}{2}\,Kb\overline \zeta \partial_t\bigg),\qquad -a\,\mathrm{Im}(\zeta\partial_\zeta) -\frac{1}{2}\,Kb\partial_t ,\\ &\mathrm{Im}\bigg(\partial_\zeta + \frac{a}{2}\,\overline{\zeta}^2\partial_{\overline\zeta} - \frac{i}{2}\,Kb\overline \zeta \partial_t\bigg),\qquad \phantom{-a\,} \mathrm{Im}(\zeta\partial_\zeta).
		\end{split}
		\end{equation*}
	These are generators of $\mathfrak{o}(2)\ltimes\mathfrak{heis}_3(\mathbb R)$ if $a=0$, $\mathfrak{u}(2)$ if $a>0$, and $\mathfrak{u}(1,1)$ if $a<0$, in complete agreement with our conclusions from the case-by-case analysis.
\end{myproof}

The attentive reader may also have noticed that the examples with positive scalar curvature are incomplete, and well-known spaces like $S^4$ or $\CP^2$ are recovered only upon adding points. This is no accident: A quaternionic K\"ahler manifold of positive scalar curvature that arises from the HK/QK correspondence is necessarily incomplete. Indeed, it is known that no complete quaternionic K\"ahler manifold of positive scalar curvature carries an almost complex structure compatible with the quaternionic structure \cite{AMP1998} or in fact, with the exception of the Grassmannians of complex $2$-planes, any almost complex structure at all \cite{GMS2011}. \Cref{prop:cpxstr}, on the other hand, shows that any example that arises from the HK/QK correspondence admits an (integrable) complex structure compatible with the quaternionic structure.

Using the expression \eqref{eq:curv-norm} for the curvature norm, we can rule out the possibility of finding a smooth completion in many cases:

\begin{prop}\label{prop:cohom}
		The metric $g^{a,b,c}$ of \eqref{eq:Liouville-QK} cannot be extended to a complete metric whenever $bc(b^2-4ac)\neq 0$ and $b\rho + 2c \le 0$ anywhere in the closure of the allowed range of the coordinate $\rho$.
\end{prop}
\begin{myproof}
	First of all, by the previous proposition, the condition $bc(b^2-4ac)\neq 0$ implies that the  metric $g^{a,b,c}$ is not locally symmetric.
	
	If there is a solution to $b\rho + 2c < 0$ in the closure of the allowed range of the coordinate $\rho$, then, because the condition is open, there is a solution  $b\rho + 2c < 0$ in the allowed range of the coordinate $\rho$. As the sign of the scalar curvature is the opposite of the sign of $b\rho + 2c$, the scalar curvature is positive. It is well-known that there are no complete quaternionic K\"ahler metrics of dimension $4$ of positive scalar curvature which are not locally symmetric \cite{FK1982, Hit1981}. Thus, $g^{a,b,c}$ cannot be extended to a complete metric.
	
	On the other hand, if there is a solution to $b\rho + 2c = 0$ in the closure of the allowed range of the coordinate $\rho$, we can find a $\rho_0$ in the neighborhood of $-\frac{2c}{b}$ such that the half-open interval between $\rho_0$ and $-\frac{2c}{b}$ (with $\rho_0$ included) is contained in the allowed range of the coordinate $\rho$.  At $\rho = -\frac{2c}{b}$, the polynomial $a\rho^2 + b\rho + c$ takes the value $\frac{c}{b^2}(4ac -b^2)\neq 0$. As a result, we have
	\begin{equation*}
		\int_{\rho_0}^{-\frac{2c}{b}}\sqrt{\frac{b\rho + 2c}{a\rho^2 + b\rho + c}}\frac{\d \rho}{\rho} < \infty.
	\end{equation*}
	By \eqref{eq:curv-norm}, $\rho = -\frac{2c}{b}$ is a curvature singularity. It follows from the finiteness of the above integral that the curvature singularity is a finite distance away from $\rho =\rho_0$. So, $g^{a,b,c}$ cannot be extended to a complete metric in this case either.
\end{myproof}

\end{document}